\documentclass[12pt,a4paper]{article}
\usepackage[utf8]{inputenc}
\usepackage{amsmath}
\usepackage{amsfonts}
\usepackage{amssymb}
\usepackage{graphicx}
\usepackage{mathabx}

    \usepackage{pdfsync}

\newtheorem{theorem}{Theorem}[section]

\newtheorem{corollary}{Corollary}[section]
\newtheorem{proposition}{Proposition}[section]

\author{Yves Le Jan}
\title{On Markovian random networks } 
\begin{document}
\maketitle

\begin{abstract}
We investigate random Eulerian networks defined by Markov loops and the associated fields,  flows and maps.\\

\end{abstract}

\section{Introduction}

The main purpose of this paper is to show that random networks, which can be defined as images of Poissonian loop ensembles, are naturally associated to random configurations and maps. In addition we present a few additional results about these networks and the flows they define.\\
% The main new results are theorems 3-1, 3-2, 3-3 and 4-1. \\
We first present briefly the framework of our study,  described in \cite{stfl}. Consider a graph $\mathcal{G}$, i.e. a set of vertices $X$ together with a set of non oriented edges $E$. We assume it is connected, and that there is no loop-edges nor multiple edges. The set of oriented edges, denoted $E^{o}$, is viewed as a subset of $X^{2}$.
 An oriented edge $(x,y)$\ is defined by the choice of an
ordering in an edge $\{x,y\}$. 
%We set $-(x,y)=(y,x)$ and if $e=(x,y)$, we denote it also $(e^{-},e^{+})$.
% The degree $d_{x}$ of a vertex $x$ is by definition the number of non oriented edges incident at $x$.

 Given a graph $\mathcal{G}=(X,E)$, a set of non negative conductances $C_{x,y}=C_{y,x}$ indexed by the set of edges $E$ and a non negative killing measure $\kappa$ on the set of vertices $X$, we can associate to them an energy (or Dirichlet form) $\mathcal{E}$, we will assume to be positive definite, which is a transience assumption. For any function $f$ on $X$, we have: $$\mathcal{E}(f,f)=\frac{1}{2}\sum_{x,y}C_{x,y}(f(x)-f(y))^{2}+\sum_x \kappa_x f(x)^{2}.$$
There is a duality measure $\lambda$ defined by $\lambda_x=\sum_y C_{x,y} +\kappa_x$.
Let $G^{x,y}$  denote the symmetric Green's function associated with $\mathcal{E}$. Its inverse equals $M_{\lambda}-C$ with $M_{\lambda}$ denoting the diagonal matrix defined by $\lambda$.\\
% It is assumed that $\sum_x G_{x,x}\lambda_x$ is finite.\\
The  associated symmetric continuous time Markov process can be obtained from the Markov chain defined by the transition matrix $P^x_{y}=\frac{C_{x,y}}{\lambda_y}$ by adding independent exponential holding times of mean $1$ before each jump. If $P$ is submarkovian, the chain is absorbed at a cemetery point $\Delta$. If $X$ is finite, the transition matrix is necessarily submarkovian.\\
%The Dirichlet form and the associated Markov process can be naturally lifted to any non ramified covering.
We denote by $\mu$ the loop measure associated with this symmetric Markov process. %It can also be viewed as a shift invariant measure on based loops 
(see \cite{stfl}).
% for the general definition in terms of Markovian bridges, but let us mention that:\\
%- the measure of a non-trivial discrete loop is the product of the transition probabilities of its edges if it is aperiodic; otherwise this product should be divided by the multiplicity of the loop.\\
%- the measure on continuous time loops is then obtained by including exponential holding times, except for one point loops on which the holding time measure (which has infinite mass)  has density $\frac{e^{-t}}{t}$.

 The Poissonian loop ensemble $\mathcal{L}_{\alpha}$ is the Poisson process of loops of intensity $\alpha\mu$, constructed in such a way that the the set of loops  $\mathcal{L}_\alpha$ increases with $\alpha$. We will denote by $\mathcal{L}_{dis,\alpha}$ the associated set of discrete time loops, which does not include one point loops.\\
 % We set $\mathcal{L}=\mathcal{L}_1$
 % Recall that when $\mathcal{G}$ is finite, $\mathcal{L}$ can be sampled by Wilson algorithm (Cf: \cite{stfl}, \cite{chang}).\\
We denote by $\hat{\mathcal{L}}_{\alpha}$ the occupation field associated with $\mathcal{L}_\alpha$ i.e. the total time spent in $x$ by the loops of $\mathcal{L}_\alpha$, normalized by $\lambda_x$.\\
 The complex (respectively real) free field is the complex (real) Gaussian field on $X$ whose covariance function is $2G$ ($G$). We will denote it by $\varphi$ (respectively $\varphi ^{\mathbb{R}}$).\\  
It has been shown in $\cite{aop} $ (see also $\cite{stfl}$) that the fields $\hat{\mathcal{L}_1}$ and $\frac{1}{2} \varphi^2$  have the same distribution. The same property holds for $\hat{\mathcal{L}}_{\frac{1}{2}}$ and $\frac{1}{2} (\varphi^{\mathbb{R}})^2$.
% Note that this property extends naturally to symmetric Markov processes in which points are %non-polar and in particular to one dimensional diffusions (see \cite{Lupdif}). Generalisations to %dimensions 2 and 3 involve renormalization (Cf \cite{stfl}).

  % Determinantal process. Grassmann field. Supersymmetry.\\
   %In what follows, we will assume for simplicity that $\mathcal{G}$ is finite. We will now define %the edge occupation fields associated with the loop ensembles.\\
Given any oriented edge $(x,y)$ of the graph, we denote
%denote by $N_{x,y}(l)$ the total number of jumps made from $x$ to $y$ by the loop $l$ and
 by $N^{(\alpha)}_{x,y}$ the total number of jumps made from $x$ to $y$ by the loops of $\mathcal{L}_{\alpha}$ and we set $N^{(\alpha)}_{\{ x,y \}} =N^{(\alpha}_{x,y}+N^{(\alpha}_{y,x}.$
 Recall that $E(\hat{\mathcal{L}}_{\alpha})=\alpha G^{x,x}$ and $E(N^{(\alpha)}_{x,y})=\alpha C_{x,y}G^{x,y}$.\\
 Recall finally the relation between $\varphi$  and the pair of fields $(\hat{\mathcal{L}_1}, (N^{(1)})$ given in Remark 12  of \cite{stfl}, Chapter 6, and the relation between $\varphi^{\mathbb{R}}$  and the pair of fields $(\hat{\mathcal{L}_1}, (N_{\{  \} } ^{(\frac{1}{2})})$ given in Remark 11:\\
  For any complex matrix $s_{x,y}$ with finitely many non zero 
 entries, all of modulus less or equal to $1$, and any finitely supported positive measure $\chi$,
 \begin{equation}\label{phi}
 E(\prod_{x,y} s_{x,y}^{N^{(1)}_{x,y}} \prod_x e^{-\sum_{x} \chi_{x} \hat{\mathcal{L}}_1^x})=E(e^{-\frac{1}{2}\sum_{x,y} C_{x,y}( s_{x,y}-1)\varphi_x \bar \varphi_y-\frac{1}{2}\sum \chi_{x}\phi_x \bar \phi_x}).
 \end{equation} 
 
% $$E(\prod_{x\neq y} Z_{x,y}^{N_{x,y}}\prod_{x} Z_{x,x}^{-(N_{x}+1)}) =E(e^{\sum_{x\neq y}(\frac{1}{2} C_{x,y} (Z_{x,y}-1)\varphi_x \bar{\varphi}_y)+\sum_{x}(\frac{1}{2} \lambda_x (1-Z_{x,x})\varphi_x \bar{\varphi}_x)} ).$$
 For any real symmetric matrix $s_{x,y}$ with finitely many non zero 
 entries, all in $[0,1)$, and any finitely supported positive measure $\chi$,
 \begin{equation} \label{phiR}
 E(\prod_{x,y} s_{x,y}^{N^{(\frac{1}{2})}_{x,y}} \prod_x e^{-\sum_{x} \chi_{x} \hat{\mathcal{L}}_{\frac{1}{2}}^x})=E(e^{-\frac{1}{2}\sum_{x,y} C_{x,y}( s_{x,y}-1)\varphi^{\mathbb{R}}_x  \varphi^{\mathbb{R}}_y-\frac{1}{2}\sum \chi_{x}({\varphi^{\mathbb{R}}_x})^2}).
 \end{equation}

  \textbf{Remarks:}\\
   - A consequence of \eqref{phi} is that for any set $(x_i,y_i)$ of distinct oriented edges, and any set $z_l$  of distinct vertices,
% \textit{distinct} from the $(x_i$ and $y_i)$'s,  
 \begin{equation}\label{phi2}
E(\prod_i N^{(1)}_{(x_i,y_i)}\prod_l (N^{(1)}_{z_l}+1))=E(\prod_i \dfrac{1}{2}C_{(x_i,y_i)}{\varphi_{x_i }\bar{\varphi}_{y_i}} \prod_l \frac{1}{2}{\lambda_{z_l } \varphi_{z_l }\bar{\varphi}_{z_l}} ).
\end{equation} 

%=\text{Per}(G_{(x_i,y_j)}){\prod_i C_{(x_i,y_i)}}{\prod_l \lambda_{z_l } G_{(z_l,z_l)}}.$$
- In particular, if $X$ is assumed to be finite, if $[D_{N^{(1)}}]_{(x,y)}=0$ for $x\neq y$ and $[D_{N^{(1)}}]_{(x,x)}=1+N^{(1)}_{x}$, for all $\chi\geq \lambda$, then (cf.\cite{lejanito}):
$$E(\det( M_{\chi}D_{N^{(1)}}-N^{(1)}))=2^{-\vert X\vert}E(\det(M_{\varphi}(M_\chi-C)M_{\bar{\varphi}})) =\det(M_\chi-C) \text{Per}(G).$$\\
 - A consequence of \eqref{phiR} is that for any set $\{x_i,y_i\}$ of distinct edges, and any set $z_l$  of distinct vertices,
% \textit{distinct} from the $(x_i$ and $y_i)$'s, 
 \begin{equation}\label{phiR2} 
E(\prod_i N^{(\frac{1}{2})}_{\{x_i,y_i\}}\prod_l (N^{(\frac{1}{2})}_{z_l}+1))=E(\prod_i C_{(x_i,y_i)}{\varphi^{\mathbb{R}}_{x_i }\varphi^{\mathbb{R}}_{y_i}} \prod_l \frac{1}{2} \lambda_{z_l } ({\varphi^{\mathbb{R}}_{z_l }})^2 ).
 \end{equation} 
%=\text{Per}(G_{(x_i,y_j)}){\prod_i C_{(x_i,y_i)}}{\prod_l \lambda_{z_l } G_{(z_l,z_l)}}.$$
- In particular, if $X$ is assumed to be finite, if $[D_{N^{(\frac{1}{2})}}]_{(x,y)}=0$ for $x\neq y$ and $[D_{N}^{(\frac{1}{2})}]_{(x,x)}=1+N^{(\frac{1}{2})}_{x}$, for all $\chi\geq \lambda$, then
$$E(\det(2 M_{\chi}D_{N^{(\frac{1}{2})}} -N_{\{  \} } ^{(\frac{1}{2})}))=E(\det(M_{\varphi^{\mathbb{R}}}(M_\chi-C)M_{\varphi^{\mathbb{R}}})) =\det(M_\chi-C) \text{Per}(G).$$
Similar expressions can be given for the minors.\\

 - Using (for example) equation (\ref{phi2}) and  (\ref{phiR2}), we can study correlation decays. For example, if $(x,y)$ and $(u,v)$ are non adjacent edges, $$Cov( N^{(1)}_{x,y},N^{(1)}_{u,v})=E(:\varphi_x \bar \varphi_y:\; :\varphi_u \bar \varphi_v:)=G^{x,v}G^{y,u}$$

- Note that a natural coupling of the free field with the loop ensemble of intensity $\frac{1}{2}\mu$ has been given by T. Lupu \cite{Lup}, using the vertex occupation field and the partition of $X$ defined by the zeros of the edge occupation field.\\

 - Most results are proved for finite graphs with non zero killing measure. Their extension to infinite transient graphs is done by considering the restriction of the energy to functions vanishing outside a finite set $D$, i.e. the Markov chain killed at the exit of $D$ and letting $D$ increase to $X$, so that the Green function of the restriction, denoted $G_D$, increases to $G$.
  %3Note that $N^{(\alpha)}_{x,x}=0$.
%\\ Let $Z$ be any Hermitian matrix indexed by pairs of vertices and $\chi$ a non-negative measure on $X$.\\
%The content of the following lemma appeared already in chapter 5 and 6 of \cite{stfl} (see remarks 11 and 13 for ii) and iii)).
%\begin{lemma}\label{toto} 
%Denote by $P^Z_{x,y}$ the matrix $P_{x,y}Z_{x,y}$.\\
%\begin{enumerate}
%\item[i)]  We have: $$E(\prod_{x\neq y} Z_{x,y}^{N^{(\alpha)}_{x,y}}e^{-\sum_x \chi_x \hat{\mathcal{L}}_{\alpha}^x} )=\left[\frac{\det(I-\frac{\lambda}{\lambda+\chi}P^{Z})}{\det(I-P)}\right]^{-\alpha}.$$
%\item[ii)] For $\alpha=1$,  $$E(\prod_{x\neq y} Z_{x,y}^{N^{(1)}_{x,y}}e^{-\sum_x \chi_x \hat{\mathcal{L}}_{1}^x} ) =E(e^{\sum_{x\neq y}(\frac{1}{2} C_{x,y} (Z_{x,y}-1)\varphi_x \bar{\varphi}_y)}e^{-\frac{1}{2}\sum_x \chi_x\varphi_x\bar{\varphi}_x} ).$$
%\item[iii)] For $\alpha=\frac{1}{2}$,  $$E(\prod_{x\neq y} Z_{x,y}^{N^{(\frac{1}{2})}_{x,y}}e^{-\sum_x \chi_x \hat{\mathcal{L}}_{\frac{1}{2}}^x} ) =E(e^{\sum_{x\neq y}\frac{1}{2} C_{x,y} (Z_{x,y}-1)\varphi ^{\mathbb{R}}_x\varphi ^{\mathbb{R}}_y}e^{-\frac{1}{2} \sum_{x}\chi_x (\varphi ^{\mathbb{R}}_x)^2} ).$$
%\end{enumerate}
%\end{lemma} 

\section{Eulerian networks}
We define a network to be a  $\mathbb{N}$-valued function defined on oriented edges of the graph. It is given by a matrix $k$ with  $\mathbb{N}$-valued coefficients which vanishes on the diagonal and on entries $(x,y)$ such that $\{x,y\}$ is not an edge of the graph. We say that $k$ is Eulerian if $$ \sum_y k_{x,y}= \sum_y k_{y,x}.$$ For any Eulerian network $k$, we define $k_x$ to be $\sum_y k_{x,y}=\sum_y k_{y,x}$.
It is obvious that the field $N^{(\alpha)}$ defines a random network which verifies the Eulerian property.\\
%Note that for a finite graph, given $N^{(\alpha)}=k$, all $\frac{\prod_x {k_x}!}{\prod_{x,y} k_{x,y}!} $ discrete time loops configurations are equally likely. Similarly, all continuous time loop configurations are also equidistributed given the values of both occupations fields.\\

 Note that $\sum_y  N^{(\alpha}_{\{ x,y \}}$ is always even.
We call even networks the sets of numbers attached to non oriented edges such that $k_x=\dfrac{1}{2}\sum_y k_{\{ x,y \}}$ is an integer.
%The distribution of the random network defined by $\mathcal{L}_{\alpha}$  was given in \cite{lejanito}. 

The cases $\alpha=1$ and $\alpha=\dfrac{1}{2}$ are of special interest. We need to recall the results of \cite{lejanito} and \cite{gaugeloops} which can be deduced from equations (\ref{phi}) and (\ref{phiR}):
\begin{theorem}\label{net}

i) For any Eulerian network $k$,$$  P(N^{(1)}=k)=\det(I-P)\frac{\prod_x {k_x}!}{\prod_{x,y} k_{x,y}!} \prod_{x,y} P_{x,y}^{k_{x,y}}.$$
ii) For any Eulerian network $k$, and any nonnegative function $\rho$ on $X$  $$P(N^{(1)}=k\:,\: \hat{\mathcal{L}}_{1}\in (\rho, \rho+d\rho ))=\det(I-P)\prod_{x,y} \frac{(\sqrt{\rho_x} C_{x,y}\sqrt{\rho_y})^{k_{x,y}}}{ k_{x,y}!}\prod_{x} \lambda_x e^{-\lambda_x\rho_x }d\rho_x.$$
 iii) For any even network $k$,$$  P(N^{(\frac{1}{2})}_{\{  \} } =k)=\sqrt{\det(I-P)}\frac{\prod_x {2 k_x}!}{\prod_x  2^{k_x} {k_x}! \prod_{x,y} k_{\{ x,y \}}!} \prod_{x,y} P_{x,y}^{k_{x,y}}.$$
 iv)  For any even network $k$, and any nonnegative function $\rho$ on $X$  $$P(N^{(\frac{1}{2})}_{\{  \} } =k\:,\: \hat{\mathcal{L}}_{\frac{1}{2}}\in (\rho, \rho+d\rho ))=\sqrt{\det(I-P)}\prod_{x,y} \frac{(\sqrt{\rho_x} C_{x,y}\sqrt{\rho_y})^{k_{\{x,y\}}}}{ k_{\{x,y\}}!}\prod_{x}\frac{\sqrt{\lambda_x}}{\sqrt{2\pi \rho_x}}e^{-\frac{1}{2}\lambda_x\rho_x }d\rho_x.$$

\end{theorem}

\textbf{Remarks:}\\

- It follows from iv) that the symmetrized $N^{(\frac{1}{2})}$ field conditioned by the vertex occupation field is, as it was observed by Werner
in \cite{wernersemiprob}, a random current model. \\

- These results can be extended to infinite transient graphs as follows: For i) given any finitely supported Eulerian network $k$ and any bounded function $F$ on Eulerian networks: 
$$E(F(N^{(1)}+k))=E(1_{\{ N^{(1)}\geq k \}}F(N^{(1)})\prod_x \frac{(N^{(1)}_x-k_x)!}{N^{(1)}_x !}
\prod_{x,y}\frac{N^{(1)}_{x,y}!} {(N^{(1)}_{x,y}-k_{x,y})!}\prod_{x,y} )P_{x,y}^{-k_{x,y}}.$$\\

- This quasi invariance property can be used to prove the closability and express the generator of some Dirichlet forms defined naturally on the space $\frak{E}$ of Eulerian networks: If $\nu$ is a bounded measure and $G$ a bounded function on $\frak{E}$, the energy of $G$ can be defined as $\int E((G( (N^{(1)}+k)-G(N^{(1)}))^2)\nu(dk)$. They define stationnary processes on $\frak{E}$, with invariant distribution given by the distribution of $N^{(1)}$.\\
Similar quasi invariance properties and stationnary processes can be derived from ii), iii) and iv). \\

\textbf{Markov property}

From theorem \ref{net}  follows a Markov property which generalizes the reflection positivity property proved in chapter 9 of \cite{stfl}:
\begin{theorem}
 Let $X$ be the disjoint union of  $X_i,\; i=1,2$ and  ${\cal{G}}_i$ be the restriction of $\cal{G}$ to $X_i$. 
 
  i) Given the values of $N^{(1)}_{x,y}$ and $N^{(1)}_{y,x}$  for $x\in X_1$ and  $y\in X_2$, the restrictions of the fields $(N^{(1)}, \widehat{\cal{ L}}_1)$ to ${\cal {G}}_1$ and ${\cal {G}}_2$ are independent.
  
 ii) Given the values of $N_{\{x,y\}}^{(\frac{1}{2})}$   for $x\in X_1$ and  $y\in X_2$, the restrictions of the fields $N_{\{\}}^{(\frac{1}{2})}, \widehat{\cal{ L}}_{\frac{1}{2}})$ to ${\cal{G}}_1$ and ${\cal {G}}_2$ are independent.
 \end{theorem}

 In both cases, we can check on the expressions given in \ref{net} ii) and iv) that after fixing the values of the conditioning, the joint density function factorizes. 
See \cite{wernersemiprob} and also \cite{camialis}  in the context of non backtracking loops.\\

 In the case $\alpha=1$, the Markov property of these fields is preserved if we modify $P$ by a factor of the form $\prod_x e^{-\beta \Phi_x}$  and a normalization constant, with $\beta >0$ and $\Phi_x$ a non-negative function of $(N^{(1)}_{x,y}, N^{(1)}_{y,x}, \widehat{\cal{ L}}_1 ^x )$\\
 
 In the case $\alpha=\frac{1}{2}$, the Markov property of these fields is preserved if we modify $P$ by a factor of the form $\prod_x e^{-\beta \Phi_x}$ and a normalization constant, with $\beta >0$ and $\Phi$ a non-negative function of $ (N^{(\frac{1}{2})}_{\{x,y\}}, \widehat{\cal{ L}}_{\frac{1}{2}} ^x )$\\
 
\textbf{Remark:} 
 An example of interest, in the case $\alpha=1$ is $\Phi_x = R_x:= (N^{(1)}_{x} )^2 - \sum_y(N^{(1)}_{x,y}) ^2$.
\section{Networks and Maps}
The distribution of $N^{(1)}$ can be interpreted further when $\mathcal{G}$ is finite. We introduce a specific type of configuration model. Let $ \mathfrak{C}$ be the space of configurations $c$ defined as follows: attach to each vertex $x$ $c_x$ entering half edges and $c_x$ exiting half edges, numbered from $1$ to $c_x$. Then provide a coupling between entering and exiting half-edges in such a way that they form an oriented edge of $\mathcal{G}$. For each oriented edge $(x,y)$, there are $c_{x,y}$ half edges exiting from $x$ coupled with $c_{x,y}$ half edges entering in $y$. For each vertex $x$, $\sum_y c_{x,y}=\sum_y c_{y,x}$ is denoted $c_x$ and we see that such a configuration $c$ defines a Eulerian network $\tilde{c}$. Moreover each Eulerian network $k$ is the image by this projection map $c\longrightarrow \tilde{c}$ of $\frac{\prod_x ({k_x}!)^2}{\prod_{(x,y)\in E^o}k_{x,y}!}$ different configurations. Indeed, there are $\prod_x \frac{k_x!}{\prod_{y}k_{x,y}!}$ way of partitioning exiting half edges, $\prod_x \frac{k_x!}{\prod_{y}k_{y,x}!}$ ways of partitioning entering half edges and $\prod_{(x,y)\in E^o}k_{x,y}!$ ways to couple them into the same oriented edge of $\mathcal{G}$. Alternatively we can partition exiting half edges only and see there are $\prod_x ({k_x}!)$ way to couple them with the entering edges.\\
We say two configurations are opposite if they are exchanged by reversing the orientation of all half edges. We say they are equivalent if they can be exchanged by a circular permutation at each vertex, acting simultaneously on entering and exiting half-edges.

 We can now conclude easily that:
\begin{theorem}
The measure $Q$ defined on $ \mathfrak{C}$ by:
$$Q(c)=\det(I-P)\frac{1}{\prod_x {c_x}!} \prod_{x,y} P_{x,y}^{c_{x,y}}$$
is a probability which projects on the distribution of the random Eulerian network $N^{(1)}$.\\
\end{theorem}

Taking unit conductances and variable killing rates, we deduce from the fact that $Q$ is a probability an expression for counting configurations with given vertex degrees: 
\begin{corollary}Denote $ \mathfrak{C}_{i_x,x\in X}$ the set of configuration with given vertex degrees $i_x$:   $ \mathfrak{C}_{i_x,x\in X}=\{c \in \mathfrak{C}, \forall x\in X, \;c_x=i_x\}$
Then the joint exponential generating function of their cardinalities $\vert \mathfrak{C}_{i_,x\in X}\vert$ is given by the inverse of the determinant of the matrix $\delta_{x,y}-s_x A_{x,y}$, A being the adjacency matrix of $\cal{G}$. 
\end{corollary}

\textbf{Configurations and Wilson algorithm:}\\
Let us call exit (entrance) configuration a class of configurations with the same partition of exiting (entering) half edges. The image of $Q$ on the set of exit configurations $ \mathfrak{C}^{ex}$ is  $Q(c^{ex})=\det(I-P) \prod_{x,y} P_{x,y}^{c_{x,y}}$. An entrance configuration is the opposite of an exit configuration. 
\\ Such an exit configuration $c^{ex}$, with the choice of an order on $X$, allows to construct a family of loops $l_{x_{i}}$ based at $x_{i}$, which do not visit $\{x_1, ..., x_{i-1}\}$, the first loop being based at the first vertex $x_1$ and visiting it  $c^{ex}_{x_{1}}$ times. It starts with the first exiting half-edge, continues to the vertex to which it is associated, then with the first half edge exiting this vertex and so on, using  and taking at each vertex exiting half edges in increasing order until all have been used at $x_1$. Then iterate this procedure with the remaining configuration, starting at the first vertex it contains, according to the order initially defined on $X$.\\
This correspondence is a bijection, and the distribution obtained on the sequence of based loops is the same as in the extension of Wilson algorithm defined in section 8-2 of \cite{stfl}. It is uniform conditionally on the network.\\
 These based loops $l_{x_{i}}$ can be randomly divided at their base points $x_i$ as in the extension of Wilson algorithm. This partition is defined by a Poisson-Dirichlet distribution if we consider loops indexed by continuous time, with exponential holding time (see remark 21 in section 8-2 of \cite{stfl}). In \cite{topo} we mentioned the existence of a discrete version of this splitting method:  If $l_{x_{i}}$ visits $x_i$ $n_{x_{i}}$ times, we can partition randomly this set of excursions  according to the exchangeable partition probability function (EPPF) $\frac{\prod_1^k (n_i -1)!}{n_{x_{i}}!}$ occurring in the so-called Chinese restaurant process (\cite{P}). Note that $\frac{\prod_1^k (n_i -1)!}{n_{x_{i}}!}$ is the probability of each partition into a set of size $n_1$, a set of size $n_2$ etc, regardless of the order. The total probability of such partitions, in the given order, is $\frac{1}{k! \prod_1^k n_i }$. 
 \\Then we can define $\mathcal{L}_{ex}$ to be the associated set of (unbased) discrete time loops.  In this way we can get a sample of the Poissonian ensemble $\mathcal{L}_{dis,1}$. Similarly, we can define $\mathcal{L}_{in}$ and a configuration $c$ determines therefore two sets of loops which are clearly independent conditionally to the network $ \tilde{c}$.\\ 
 
 \textbf{Configurations and maps:}\\
  The relations between configuration and Poissonian loop ensembles appears to be deeper, as shown in the following.\\
 - Recall that a combinatorial map can be defined as a graph, with possibly multiple edges between vertices, equipped with a combinatorial imbedding, i.e. a cyclic ordering of edges around each vertex (see \cite{maps}, \cite{LanZvo}). This ordering allows to define faces and the map can be drawn on a surface whose genus (more precisely the sum of the genera of its connected components) is canonically defined in terms of the Euler-Poincar\'e characteristic. Let us say it is a $\cal{G}$-map if its vertex set is contained in $X$ and if its edges are multiple copies of elements of $E$ (if $\cal{G}$ is complete, the second condition is redundant). Let us say a map is numbered if a first edge is chosen at every vertex.\\
 - An element $c$ of $\mathfrak{C}$  defines a numbered $\cal{G}$-map $\mathcal{M}(c)$ with an even number of incident edges at each vertex. The first edge at each vertex is determined by the first exiting half edge and the half edge it is coupled to, the second by the first entering half edge and the  half-edge it is coupled to, the third by the second exiting half edge and the half-edge it is coupled to,  and so on, alternating entering and exiting half edges according to their original cyclic order. We get in this way, at each vertex, a cyclic order on edges, alternating exiting and entering half-edges.
 Two equivalent configurations define the same $\cal{G}$-map.\\ 
 -We now show that an equivalent class of configurations defines also a pair $(\mathcal{L}_{+},\mathcal{L}_{-})$ of samples of $\mathcal{L}_{dis,1}$ and that we can define alternating signs on the faces of the map so that these sets of loops are respectively the positive (negative) faces of the map.\\
  $\mathcal{L}_{+}$ is defined by the cycles of the permutation on exiting half-edges obtained by mapping an exiting half-edge with the exiting half-edge following, (in the above defined cyclic order at each vertex), the entering half-edge coupled to it. It induces the network $\tilde{c}$\\
  $\mathcal{L}_{-}$ is defined by the cycles of the permutation on all entering half-edges obtained by mapping an entering half-edge with the entering half-edge following the exiting half-edge coupled to it. As it uses the coupling backwards, it induces the opposite network $-\tilde{c}$\\
Each edge of the map is in this way adjacent to two faces, the positive one defined by exiting half edges, and the negative one by entering half edges. Note that $\mathcal{L}_{+}$ determines the coupling, hence the configuration. The same holds for $\mathcal{L}_{-}$. \\ 
 - We can now state the following :
  \begin{theorem}The image on $\mathcal{G}$ of the sets of oriented face contours $\mathcal{L}_{+}$ and $\mathcal{L}_{-}$ have both the same distribution as $\mathcal{L}_{dis,1}$ .
  \end{theorem}
  We can again choose an arbitrary order on the vertices and use the coupling and the ordering at each vertex to run a different version of Wilson algorithm to construct $\mathcal{L}_{+}$. This algorithm can be viewed as a way, starting from a given exit configuration (i.e. a partition of exiting half edges) to sample the random coupling with entering edges progressively. We start a based loop at the first point with the first exiting half edge, then move to the vertex attached to it by the partition, and choose an entering half edge to couple it. Then follow the exiting half edge following this entering edge in the alternating cyclic order. At every return to the base point, the based loop is created if the first exiting half-edge follows the last entering one, in the cyclic order. Conditioning on the initially given partition of exiting half edges and on the partial coupling with entering half edges done until this return time, this occurs with probability one over the number of unused entering half edges at the base point. Then restart at the same vertex with the next free exiting half-edge. When none is left, restart the algorithm at the next possible vertex with the first unused exiting half-edge. The concatenation of the based loops obtained at each vertex has the same distribution as in Wilson algorithm.The random ordered partitions of these loops we have also obtained (in contrast with the previous version of the algorithm) is exactly the size-biased version of the ones defined by the above mentioned EPPF $\frac{\prod_1^k (n_i -1)!}{n_{x_{i}}!}$ (see section 2-1 in \cite{P}). This observation completes the proof of  the theorem as the case of $\mathcal{L}_{-}$ can be treated in the same way.\\

\textbf{Euler-Poincar\'e characteristic}\\
The Euler-Poincar\'e characteristic $\chi (c)$ of the map defined by a configuration $c$ is $ \vert \{ x\in X, \, c_x>0\} \vert + \vert \mathcal{L}_{+}(c)\vert  + \vert \mathcal{L}_{-}(c)\vert  - N(c)$ with $N(c)=\sum_{x,y} \tilde{c}_{x,y}$. It is the sum of the Euler-Poincar\'e characteristics of its connected components.
%\begin{proposition} The distribution of $\chi$ is determined by its  generating function:
%$$E(u^{\chi})=\sum_{k\in \frak{E}} u^{-\sum_x (k_x-1)^+}\frac{E(u^{\vert \mathcal{L}_{dis,1}\vert }1_{N^{(1)}=k})^2}{P(N^{(1)}=k)}$$
%\end{proposition}
%
%This is a straightforward consequence of the definition of $\chi $ and of the previous theorem.
%Note that  $E(u^{\vert \mathcal{L}_{dis,1}\vert }1_{N^{(1)}=k})$ can be computed as $\det(I-P)^{1-u}P(N^{(u)}=k)$. Indeed, $E(u^{\vert \mathcal{L}_{dis,1}\vert }\prod_{x,y}s_{x,y}^{N_{x,y}^{(1)}})=e^{\sum_l( u \prod _{x,y}s_{x,y}^{N_{x,y}(l)}-1) \mu(l) }=\det(I-P)^{1-u}E(\prod_{x,y}s_{x,y}^{N_{x,y}^{(u)}}) $\\
%
%\medskip
%  \textbf{Remarks:}\\
%  
%  
%- Recall from \cite{lejanito} that $ \det(I-P)^{-u}P(N^{(u)}=k)=C \prod_{x,y} P_{x,y}^{k_{x,y}}$ where $C$ is the coefficient of $\prod_{x,y} P_{x,y}^{k_{x,y}} $ in the $u$-permanent $\mathrm{Per}_{u}(P(k_x, x\in X)) $. \\ Here $P(k_x, x\in X)$ is denoting the $(\vert k \vert , \vert k \vert)$ matrix obtained by repeating $k_x$ times each column of $P$ of index $x$  and then $k_y$ times each line of index $y$  (with  $\vert k \vert = \sum_{x} k_{x}$).\\
%%The $u$-permanent $\mathrm{Per}_{u}$ of a $(d,d)$-matrix $A$ is defined as the sum, on all permutations $\sigma$ of $\{1,...,d\}$, of the products $u^{\#\mathrm{cycles}(\sigma)}\prod_{i=1}^{d} A_{i,\sigma(i)}$.\\

\medskip  
 Using the results of  \cite{stfl}, the expectation of $\chi$ can be expressed as$$E(\chi)=\sum_x (1-\frac{1}{\lambda_x G^{x,x}})-2 \ln\det(I-P) -\sum_{x,y}C_{x,y}G^{x,y}.$$ Expressions can also be given for the number and the sizes of the connected components of the map (see \cite{ljlm}).\\

   All these results apply in particular to the case of complete graphs, which is of special interest in relation with the theory of combinatorial maps.\\ For $d$ vertices and constant killing rate $\kappa$, $\lambda=d-1+\kappa$ and $G=\frac{1}{d+\kappa} (I+\frac{1}{\kappa}J)$, $J$ denoting the $(d,d)$ matrix with all entries equal to $1$. Setting $p=\frac{1}{d-1+\kappa}$,  $\det( \delta_{x,y}-p A_{x,y})=(1-(d-1)p)(1+p)^{d-1}$. In particular $E(\chi)=d(1-(1-\frac{1}{\kappa+1})(1+\frac{1}{d-1+\kappa}))-\frac{d(d-1)}{\kappa(d+\kappa)}-2\ln(\frac{\kappa}{d-1+\kappa})-2(d-1)\ln(1+\frac{1}{d-1+\kappa}).$\\
   If we let $d$ and $\kappa$ increase to infinity, with $\frac{\kappa}{d}$ converging to zero, the expected number of vertices and the expected number of edges are both equivalent to  $\frac{d}{\kappa}$. More precisely, they are respectively equal to  $\frac{d}{(\kappa+1)}+1+o(1) $ and $\frac{d}{\kappa}-1+o(1) $. The expected number of faces is $2\ln(d/\kappa) -2+o(1)$. The expectation of the characteristic $\chi$ can converge to any value $v$ by taking $\kappa(d,v)=\sqrt{\dfrac{d}{u(d,v)}}$ with $u(d,v)>1$ defined as the unique solution of the equation $u(d,v)-\ln(u(d,v))=\ln(d)-v$. Note that $u(d,v)$ is equivalent to $\ln(d)$ as $d$ increases to infinity. If we take $\kappa(d)=\sqrt{\dfrac{d}{\ln(d)}}$, the expectation of $\chi$ is equivalent to $\ln(\ln(d))$ as $d$ increases to infinity. \\
   Consider now essential vertices, i.e. such that $N_x>1$. Their expected number in general is $\sum_x  (1-\frac{1}{\lambda_x G^{x,x}})^2$ which equals $\frac{d(d+2\kappa-1)^2}{(\kappa+1)^2(d-1+\kappa)^2}$ in the complete graph case. If we take  $\kappa(d)=\sqrt{\dfrac{d}{\ln(d)}}$ or $\kappa(d,v)=\sqrt{\dfrac{d}{u(d,v)}}$, the expected number of faces and the expected number of essential vertices are both equivalent to $\ln(d)$ as $d$ increases to infinity.\\
   
     %They can be related to the combinatorial theory of maps (or imbedded multigraphs): see \cite{maps}, \cite{LanZvo}.\\
 
 \textbf{Remark:}\\   
  If the graph is bipartite, namely if there is a partition of $\{X_1,X_2\}$ of $X$ such that $E^o \subseteq X_1 \times X_2 \cup X_2 \times X_1$, $\mathfrak{C}$ is the set of pairs of numbered $\cal{G}$-maps with the same number of incident edges at each vertex.\ The first map $M_1$ is defined by oriented edges entering in $X_1$ and exiting from $X_2$ and the second map $M_2$ is by oriented edges entering in $X_2$ and exiting from $X_1$.\\

\textbf{Configurations and even networks:}\\
A similar construction can be made with even networks. The space of configurations is now the set $ \mathfrak{C}^{ev}$  of numbered  $\cal{G}$-maps with an even number of incident edges at each vertex. At each vertex $x$ there are $ \frac{2k_x!}{\prod_{y}k_{\{x,y\}}!}$ way of partitioning the $2k_x$ incident edges and there are $\prod_{\{x,y\}\in E}k_{\{x,y\}}!$ ways to couple the numbers assigned to the edges between $x$ and $y$.  Each even network $k$ is the image  of $\frac{\prod_x ({2k_x}!)}{\prod_{\{x,y\}\in E}k_{\{x,y\}}!}$ different configurations. We have :
\begin{theorem}
The measure $Q^{ev}$ defined on $ \mathfrak{C}^{ev}$ by
$$Q(c)=\sqrt{\det(I-P)}\frac{1}{\prod_x {2^{c_x}c_x}!} \prod_{x,y} P_{x,y}^{c_{x,y}}$$
is a probability which projects on the distribution of the edge occupation field $(N^{(\frac{1}{2})}_e, e\in E)$.
\end{theorem}
A similar corollary holds:
\begin{corollary}Denote $ \mathfrak{C}^{(ev)}_{i_x,x\in X}$ the intersection of $\{c \in \mathfrak{C}^{(ev)}, c_x=i_x\}$
Then the joint exponential generating function of their cardinalities $\vert \mathfrak{C}^{(ev)}_{i_x}\vert ,x\in X$ is given by the inverse square root of the determinant of the matrix $\delta_{x,y}-\frac{s_x}{2}  A_{x,y}$, A being the adjacency matrix of $\cal{G}$.
\end{corollary}

- For the complete graph with constant killing rate $\kappa$ and $\alpha=\frac{1}{2}$, we see that the Eulerian network  is very close of the configuration model studied extensively in chapters 7 and 10 of \cite {rem}. 

%The expected number of occupied edges is $\frac{d(d-1)}{\kappa(d+\kappa)}$, with ${d=\vert X\vert}$. Asymptotically for large $d$, the network should have one large connected component of order $d$ and other components of strictly smaller order (see also \cite{ljlm}).\\

\section{Random homology and Flows}
\
Let us assume the graph $\mathcal{G}$ is finite. We now recall a result of \cite{lejanito}. The additive semigroup of Eulerian networks is naturally mapped on the first homology group $H_1(\mathcal{G},\mathbb{Z})$ of the graph, which is  an Abelian group with $n=\vert E \vert -\vert X \vert +1$ generators. The homology class of the network $k$ is determined by the antisymmetric part $\widecheck{k}$ of the matrix $k$.\\
 We denote by $H^1(\mathcal{G},\mathbb{R})$ the space of harmonic one-forms, which in our context is the space of one-forms $\omega^{x,y}=-\omega^{y,x}$ such that $\sum_{y}C_{x,y}\omega^{x,y}=0$ for all $x\in X$ and by $H^1(\mathcal{G},\mathbb{Z})$ the space of harmonic one-forms $\omega$ such that for all discrete loops (or equivalently for all non backtracking discrete loops) $\gamma$, the holonomy $\omega(\gamma)$ is an integer. The Jacobian torus of the graph $ Jac(\mathcal{G})$ is defined to be the quotient $H^1(\mathcal{G},\mathbb{R})/ H^1(\mathcal{G},\mathbb{Z})$. The distribution of the induced random homology $\widecheck{N}^{(\alpha)}$ can be computed as a Fourier integral on the Jacobian torus for the Lebesgue measure associated with the metric defined by the conductances.\\
We have (see \cite{lejanito} ):
 \begin{proposition} 
$$ P(\widecheck{N}^{(\alpha)}=j)=\frac{1}{\vert Jac(\mathcal{G})\vert}\int_{Jac(\mathcal{G})}\left[\frac{\det(G^{(2\pi i\omega)})}{\det(G)}\right]^{\alpha}e^{-2\pi i\langle j,\omega \rangle}d\omega.$$
 \end{proposition}
 A similar integral expression can be given for the joint distribution of $\widecheck{N}^{(\alpha)}$ and the vertex occupation field.\\

We say that a Eulerian network $j$ is a flow if it defines an orientation on edges on which it does not vanish, i.e. iff $j_{x,y}\,j_{y,x}=0$ for all edges $\{x,y\}$.
It is easy to check that the measure $j_x=\sum_y j_{x,y}$ is preserved by the Markovian matrix $q$ defined as follows:
$q^x_y =\frac{j_{x,y}}{j_x}$ if $j_x >0$, $q^x_y =\delta^x_y$ if $j_x =0$.\\ We can define the stochasticity of the flow at $x$ to be $S_x=j_x^2-\sum_y j_{x,y}^2$. If it vanishes everywhere, the Markovian transition matrix is a permutation of $X$.\\
We can define the flow $j(\widecheck{k})$ associated to an element  $\widecheck{k}$ of  $H_1(\mathcal{G},\mathbb{Z})$ by $$j_{x,y}=1_{\{\widecheck{k}_{x,y}>0\}}\,\widecheck{k}_{x,y}.$$
The mapping $j$ is clearly a bijection. \\

We now show that a simple expression of this joint distribution of the flow and the vertex occupation field  can be derived from theorem \ref{net}  when $\alpha=1$.\\
Let $\frak{J}$ be the set of flows on $\mathcal{G}$. %For $h\in \frak{J}$, set $E^h={(x,y)\in E^o, h_{x,y}>0}$\\
 and $\frak K$ be the set of $\mathbb N$-valued sets of conductances on $\cal G$.
$$\{j(\widecheck{N}^{(1)})=h\}=\bigcup _{k\in \frak K}\{\bigcap_{\{x,y\}\in E}\{ N^{(1)}_{x,y}=k_{\{x,y\}}+h_{x,y},N^{(1)}_{y,x}=k_{\{x,y\}}+h_{y,x}\}\}.$$
From \ref{net} ii), it follows that:
$$P(j(\widecheck{N}^{(1)})=h ,\: \hat{\mathcal{L}}_{1}\in (\rho, \rho+d\rho ))=\sum_{k\in \frak K}\det(I-P)\prod_{x,y} \frac{(\sqrt{\rho_x} C_{x,y}\sqrt{\rho_y})^{k_{x,y}+h_{x,y}}}{ (k_{x,y}+h_{x,y})!}\prod_{x}\lambda_x e^{-\lambda_x\rho_x }d\rho_x.$$

Recall the definition of the modified Bessel function:
$$I_{\nu}(x)=\sum_{m=0}^{\infty} \frac{1}{m!\Gamma(\nu+m+1)}\left( \frac{x}{2}\right)^{2m+\nu}. $$ \\
From this follows:
\begin{theorem}
For any $h\in \frak J$, setting $h_{\{x,y\}}=\sup(h_{x,y},h_{y,x})$
$$P(j(\widecheck{N}^{(1)})=h ,\: \hat{\mathcal{L}}_{1}\in (\rho, \rho+d\rho ))=\det(I-P)\prod_{\{x,y\}} I_{h_{\{x,y\}}}(2\sqrt{\rho_x} C_{x,y}\sqrt{\rho_y})\prod_{x}\lambda_xe^{-\lambda_x\rho_x }d\rho_x.$$
\end{theorem}

\textbf{Remarks}\\

- Using again equation (\ref{phi2}) and  (\ref{phiR2}), we can study correlation decays.  $$E(\widecheck{N}^{(1)}_{x,y},\widecheck{N}^{(1)}_{u,v})=E((\varphi_x \bar \varphi_y-\varphi_y \bar \varphi_x)(\varphi_u \bar \varphi_v-\varphi_v \bar \varphi_u))=2G^{x,v}G^{y,u}-2G^{x,u}G^{y,v}.$$ 
If  $\mathcal{G}$ is the three dimensional square lattice, this correlation will decay like the fourth power of the distance.   

We can also consider the circulation of the flow around plaquettes and get similar results, with correlation decaying like the sixth power of the distance.

- Recall that we defined $R_x:= (N^{(1)}_{x} )^2 - \sum_y(N^{(1)}_{x,y}) ^2$. Note that $R=0$ implies that the flow defined by $\widecheck{N}^{(1)}$  has zero stochasticity (but the converse is not true). Hence, as $\beta$ increases to infinity, the probability modified by $\prod_x e^{-\beta R_x}$ concentrates on the set of flows of null stochasticity (as it has positive $P$-probability).

%\textbf{Acknowledgment:} Hearty thanks are due to the referee for his careful reading and useful suggestions.

\bigskip

\noindent
NYU Shanghai. 1555 Century Blvd, Pudong New District. Shanghai. China. \\
and\\
  D\'epartement de Math\'ematique. Universit\'e Paris-Sud.  Orsay, France.\\
  
   \bigskip   
\noindent 
   yves.lejan at math.u-psud.fr\\
   and\\ 
   yl57 at nyu.edu

\end{document}